\documentclass[11pt,letterpaper]{amsart}

\usepackage[letterpaper,margin=.95in]{geometry}
\usepackage[all]{xy} 
\usepackage{amsmath, amssymb, amsfonts, latexsym, mdwlist, amsthm, amscd,mathabx,braket}
\usepackage{subfig}
\usepackage{graphicx}
\usepackage{wrapfig}
\usepackage{etex}
\usepackage{tikz-cd}
\usepackage[colorlinks=true, citecolor=blue, urlcolor=blue, linkcolor=blue, pagebackref]{hyperref}

\usepackage{tikz}
\usetikzlibrary{calc,trees,positioning,arrows,chains,shapes.geometric,%
    decorations.pathreplacing,decorations.pathmorphing,shapes,%
    matrix,shapes.symbols}

\tikzset{
>=stealth',
  punktchain/.style={
    rectangle,
    rounded corners,
    draw=black, thick,
    minimum height=3em,
    text centered,
    on chain},
  line/.style={draw, thick, <-},
  eLement/.style={
    tape,
    top color=white,
    bottom color=blue!50!black!60!,
    minimum width=8em,
    draw=blue!40!black!90, very thick,
    text width=10em,
    minimum height=3.5em,
    text centered,
    on chain},
  every join/.style={->, thick,shorten >=1pt},
  decoration={brace},
  tuborg/.style={decorate},
  tubnode/.style={midway, right=2pt},
}

\usepackage{paralist}
\usepackage{enumitem} 
\setdefaultenum{(1)}{(a)}{(i)}{}
\setlist[enumerate,1]{label={\upshape(\arabic*)}}
\setlist[enumerate,2]{label={\upshape(\alph*)},ref=\theenumi\upshape(\alph*)}
\setlist[enumerate,3]{label={\upshape(\roman*)},ref=\theenumi\theenumii\upshape(\roman*)}
\usepackage[capitalize]{cleveref}
\crefname{Prop}{Proposition}{Propositions}
\crefname{Thm}{Theorem}{Theorems}
\crefname{Lem}{Lemma}{Lemmas}
\crefname{enumi}{Case}{Cases}

\def\C{\ensuremath{\mathbb{C}}}

\def\P{\ensuremath{\mathbb{P}}}

\def\R{\ensuremath{\mathbb{R}}}

\def\Aut{\mathop{\mathrm{Aut}}\nolimits}

\def\MG13{\ensuremath{{\mathcal M}_{\Gamma_1(3)}}}
\def\tildeMG13{\ensuremath{\widetilde{\mathcal M}_{\Gamma_1(3)}}}
\def\Stab{\mathop{\mathrm{Stab}}\nolimits}

\def\Db{\mathrm{D}^{\mathrm{b}}}

\newcommand\TFILTB[3]{%
\xymatrix@=1pc{
{0 = {#1}_0} \ar[rr]&&
{{#1}_1} \ar[rr]\ar[ld] &&
{{#1}_2} \ar[r]\ar[ld] &
{\cdots} \ar[r] & { {#1}_{#3-1}} \ar[rr] &&
{{#1}_{#3} = {#1}} \ar[ld]
\\
& *{{#2}_1} \ar@{.>}[ul] &&
{{#2}_2} \ar@{.>}[ul] & &&&
{{#2}_{{#3}}} \ar@{.>}[ul]
}}

\def\abs#1{\left\lvert#1\right\rvert}

\makeatletter
\newtheorem*{rep@theorem}{\rep@title}
\newcommand{\newreptheorem}[2]{%
\newenvironment{rep#1}[1]{%
 \def\rep@title{#2 \ref{##1}}%
 \begin{rep@theorem}}%
 {\end{rep@theorem}}}
\makeatother

\newtheorem{Thm}{Theorem}
\newreptheorem{Thm}{Theorem}

\newtheorem{Cor}[Thm]{Corollary}
\newreptheorem{Cor}{Corollary}

\newreptheorem{Con}{Conjecture}

\newtheorem*{theorem*}{Theorem}
\newtheorem*{lemma*}{Lemma}
\newtheorem*{proposition*}{Proposition}
\newtheorem*{conjecture*}{Conjecture}
\newtheorem*{corollary*}{Corollary}
\newtheorem*{problem*}{Problem}

\newtheorem{Thm-int}{Theorem}

\theoremstyle{definition}
\newtheorem{Def-s}[Thm]{Definition}

\newtheorem{Rem}[Thm]{Remark}

\def\C{\ensuremath{\mathbb{C}}}

\def\P{\ensuremath{\mathbb{P}}}

\def\R{\ensuremath{\mathbb{R}}}

\def\AA{\ensuremath{\mathcal A}}

\def\XX{\ensuremath{\mathcal X}}
\def\YY{\ensuremath{\mathcal Y}}

\newcommand{\ignore}[1]{}
\begin{document}

\title{Stability conditions on some families of Calabi-Yau threefolds via orbifolding}

\author{Howard Nuer}
\address{Faculty of Mathematics \\
Technion, Israel Institute of Technology}
\email{hnuer@technion.ac.il}

\keywords{Bridgeland Stability, Calabi-Yau threefolds, Mirror Symmetry}
\subjclass[2010]{Primary: 14F05; Secondary: 14J32, 18E30.}

\begin{abstract}
We prove that families of Calabi-Yau threefolds (CY3's) admit Bridgeland stability conditions when they are obtained via orbifolding from a family of CY3's admitting Bridgeland stability conditions.  In particular, we prove that the quintic mirror admits Bridgeland stability conditions.
\end{abstract}

\maketitle

The notion of stability conditions on triangulated categories was introduced by Bridgeland in \cite{Bri:07} in order to make a mathematically rigorous definition of Douglas's $\Pi$-stability for $D$-brains in string theory.
Ever since, the existence of stability conditions on projective varieties of higher dimension, most importantly on CY3's, is considered the biggest open problem in the theory of Bridgeland stability.  
While stability conditions are known to exist now on a number of smooth families of projective threefolds, a major breakthrough came with Li's construction of geometric stability conditions on quintic threefolds in \cite{Li:StabilityQuintic}, which are the first family of strict CY3's shown to admit stability conditions.
Inspired by Li's techniques, Koseki and Liu have now constructed geometric stability conditions on CY double/triple solids \cite{Kos} and CY complete intersections of type $(2,4)$ \cite{Liu}, respectively.  In this note, we show that their constructions will induce stability conditions on any families of CY3's obtained via  ``(smooth) orbifolding''.

Given a family of smooth Calabi-Yau threefolds $\phi\colon\XX\to A$, smooth orbifolding is a common method for producing a mirror family $\psi\colon\YY\to B$ of Calabi-Yau threefolds that goes back to the original example of the quintic.
For this method, one chooses a subvariety $A_0\subset A$ parametrizing a subfamily of $\XX$ with a finite group $G$ acting non-trivially and equivariantly on $\phi_0\colon\XX_0\to A_0$, where $\phi_0$ is the restriction of $\phi$ to $\XX_0=\phi^{-1}(A_0)$.  
One then takes the quotient by $G$ to obtain a family $\bar{\phi}_0\colon\bar{\XX}_0\to B=A_0/G$ of singular Calabi-Yau threefolds whose singularities we resolve crepantly and simultaneously in the family to obtain a new family of smooth CY3's $\psi\colon\YY\to\bar{\XX}_0\to B$.  
In the more general version of the orbifolding method, one may first take  a partial compactification $\overline{\phi}\colon\overline{\XX}\to\overline{A}$ of the family $\phi$ to allow for CY3's with mild singularities and then choose the subfamily $A_0$ to lie in the boundary $\overline{A}-A$ so that the general fiber of $\phi_0\colon\XX_0\to A_0$ is a singular CY3 even before taking the quotient by $G$ (which can now also be taking to be trivial in this more general version).
The last step of the method is the same with a simultaneous, crepant resolution of the singularities of $\bar{\phi}_0$.

In this short note, we prove that in the framework of smooth orbifolding, existing methods in Bridgeland stability and derived categories allow us to induce stability conditions on the CY3's produced via smooth orbifolding from those CY3's admitting stability conditions.

\begin{Thm}\label{Thm:StabilityOnMirror}
    Let $X$ be a smooth Calabi-Yau threefold admitting Bridgeland stability conditions as in \cite{BMT14a}.  Then any smooth member $Y$ of a family of CY3's obtained via smooth orbifolding admits Bridgeland stability conditions.
\end{Thm}
\begin{proof}
    By the orbifolding assumption, there is a finite group $G$ and a birational projective morphism $f\colon Y\to X/G$ such that $K_Y=f^*K_{X/G}=0$.  It follows from \cite[Theorem 1.2]{BKR01} and \cite[Theorem 1.1]{Bri02} that $\Db(Y)\cong\Db_G(X)$, where this denotes the equivariant bounded derived category.  Thus it follows from \cite[Theorem 1.1]{MMS09} that $\Stab(X)_G$, the set of $G$-invariant stability conditions on $\Db(X)$, embeds as a closed submanifold of $\Stab(\Db_G(X))=\Stab(Y)$.  Taking a stability condition $\sigma=(\AA_{\omega,B},Z_{\omega,B}$ on $X$ as constructed in \cite{BMT14a} with $\omega$ and $B$ invariant under the $G$-action, we get $\sigma$ is a $G$-invariant stability condition on $X$ which induces a stability condition on $Y$, as required.
\end{proof}

\begin{Cor}\label{Cor:StabilityOnMirrorQuintic}
    Let $X\subset\P^4$ be a smooth quintic hypersurface and let $Y$ be a smooth CY3 in any of the 43 families constructed via smooth orbifolding from quintic CY3's in \cite{Yu:McKay}.  Then $Y$ admits Bridgeland stability conditions.  In particular, the quintic mirror $Y$ admits Bridgeland stability conditions.
\end{Cor}
\begin{proof}
By \cite[Theorem 1.3]{Li:StabilityQuintic}, each smooth quintic threefold $(X,H)$ admits a continuous family of geometric stability conditions $\sigma_{\alpha,\beta,H}^{a,b}$ parametrized by $(\alpha,\beta,a,b)\in\R_{>0}\times\R\times\R_{>0}\times\R$ such that
$$\alpha^2+(\beta-\lfloor\beta\rfloor-\frac{1}{2})^2>\frac{1}{4};\quad\text{and}\quad a>\frac{\alpha^2}{6}+\frac{1}{2}\abs{b}\alpha.$$
It follows from \cref{Thm:StabilityOnMirror} that any CY3 obtained from the quintic via smooth orbifolding admits Bridgeland stability conditions as well.  
Essentially all possibilities are covered by the 43 families described in \cite[Theorem 3.2]{Yu:McKay}.
The quintic mirror CY3 is on this list and its construction via smooth orbifolding from the quintic was nicely described for the mathematical community in \cite[Section 5]{Morrsion:Quintic}.  
\end{proof}
\begin{Rem}
    It would be interesting to determine if the Bridgeland stability conditions constructed above on the quintic mirror CY3 are also geometric.  
\end{Rem}
\begin{Rem}
    By studying the argument in \cref{Cor:StabilityOnMirrorQuintic} and writing down the locus of invariant stability conditions in detail, it should be possible to explicitly see the complex moduli space of $Y$ inside the double quotient $[\Aut(\Db(X))\backslash \Stab(X)/\C]$ where it appears as the stringy K\"{a}hler moduli space of the quintic $X$ itself.  In particular, it would be interesting to determine how this subset overlaps with the neighborhoods of the stringy K\"{a}hler moduli space described in \cite{Li:StabilityQuintic,Toda:GepnerQuintic}.
\end{Rem}

\begin{Rem}
In principal, the main impediments to extending \cref{Thm:StabilityOnMirror} to apply to the general orbifolding method are constructing stability conditions on the fibers of the family $\phi_0$ of  singular CY3's and proving a more general derived McKay correspondence allowing for the original variety to be singular.
\end{Rem}
\subsection*{Acknowledgements}
The author would first like to thank Jake Solomon for conversations and questions that inspired this note.  He would also like to thank Lev Borisov, Tyler Kelly, Chunyi Li, and Emanuele Macr\`{i} for helpful conversations that increased the scope of applications of \cref{Thm:StabilityOnMirror}.

\bibliographystyle{plain}
\bibliography{Common_Biblio}
\end{document}